\documentclass{gtart_h}

\def\ifplaintex{\expandafter\ifx\csname documentclass\endcsname\relax}

\def\gtp{{\mathsurround=0pt\it $\cal G\mskip-2mu$eometry \&\ 
$\cal T\!\!$opology $\cal P\!$ublications}}  

\def\recd{{\small Received:\qua\receiveddate\ifx\reviseddate\relax
\else\qquad Revised:\qua\reviseddate\fi\par}} 

\def\lognumber#1{\def\thelognumber{#1}}
\def\volumenumber#1{\def\thevolumenumber{#1}}
\def\volumeyear#1{\def\thevolumeyear{#1}}
\def\papernumber#1{\def\thepapernumber{#1}}
\def\pagenumbers#1#2{\def\startpage{#1}\def\finishpage{#2}}
\def\published#1{\def\publishdate{#1}}

\def\received#1{\def\receiveddate{#1}}

\def\accepted#1{\def\accepteddate{#1}}

\def\asciiaddress#1{\def\theasciiaddress{#1}}
\def\asciiemail#1{\def\theasciiemail{#1}}

\let\\\par\let\thelognumber\relax\let\thevolumenumber\relax
\let\thepapernumber\relax\let\thevolumeyear\relax\let\startpage\relax
\let\finishpage\relax\let\publishdate\relax\let\receiveddate\relax
\let\reviseddate\relax\let\accepteddate\relax\let\theasciititle\relax
\let\theasciiauthors\relax\let\theasciiaddress\relax
\let\theasciiabstract\relax

\let\theasciiemail\relax

\ifplaintex
\font\logobig=cmssbx10 scaled 3836
\font\logomed=cmssbx10 scaled 2557
\else
\font\logobig=cmssbx10 scaled 4200
\font\logomed=cmssbx10 scaled 2800
\fi

\long\def\makeagttitle{   
\count0=\startpage
\agt\hfill      
\hbox to 45truept{\vbox to 0pt{\vglue -13truept{\logomed A\kern -.37em{\logobig 
T}\kern -.38em G}\vss}\hss}
\break
{\small Volume \thevolumenumber\ (\thevolumeyear)
\startpage--\finishpage\nl
Published: \publishdate}

\vglue .25truein

{\parskip=0pt\leftskip 0pt plus
1fil\def\\{\par\smallskip}{\Large\bf\thetitle}\par\medskip} \vglue
0.05truein

{\parskip=0pt\leftskip 0pt plus 1fil\def\\{\par}{\sc\theauthors}
\par\medskip}%
 
\vglue 0.03truein 

{\small\leftskip 25truept\rightskip 25truept{\bf Abstract}\stdspace\theabstract

{\bf AMS Classification}\stdspace\theprimaryclass
\ifx\thesecondaryclass\relax\else; \thesecondaryclass\fi\par
{\bf Keywords}\stdspace \thekeywords\par}\vglue 7truept

}   

\ifplaintex
\font\phead=cmsl9 scaled 950
\font\pnum=cmbx10 scaled 913
\font\pfoot=cmsl9 scaled 950
\headline{\vbox to 0pt{\vskip -4.5mm\line{\small\phead\ifnum
\count0=\startpage ISSN 1472-2739 (on-line) 1472-2747 (printed)
\hfill {\pnum\folio}\else\ifodd\count0\def\\{ }%
\ifx\theshorttitle\relax\thetitle\else\theshorttitle\fi\hfill{\pnum\folio}
\else\def\\{ and }{\pnum\folio}\hfill\ifx\theshortauthors\relax\theauthors
\else\theshortauthors\fi\fi\fi}\vss}}
\footline{\vbox to 0pt{\vglue 0mm\line{\small\pfoot\ifnum\count0=\startpage
\copyright\ \gtp\hfill\else
\agt, Volume \thevolumenumber\ (\thevolumeyear)\hfill\fi}\vss}}
\else
\font=cmsl9 scaled 1050
\font\lnum=cmbx10 
\font=cmsl9 scaled 1050
\makeatletter
\def\@oddhead{{\small\lhead\ifnum\count0=\startpage ISSN 1472-2739 
(on-line) 1472-2747 (printed)\hfill {\lnum\number\count0}\else\ifodd\count0
\def\\{ }\ifx\theshorttitle\relax \thetitle \else\theshorttitle\fi\hfill
{\lnum\number\count0}\else\def\\{ and }{\lnum\number\count0}
\hfill\ifx\theshortauthors\relax 
\theauthors\else\theshortauthors\fi\fi\fi}}\def\@evenhead{\@oddhead}
\def\@oddfoot{\small\lfoot\ifnum\count0=\startpage\copyright\ \gtp\hfill\else
\agt, Volume \thevolumenumber\ (\thevolumeyear)\hfill\fi}
\def\@evenfoot{\@oddfoot}
\makeatother
\fi
\let\maketitlepage\makeagttitle
\let\makeshorttitle\maketitlepage
\let\maketitle\maketitlepage

\newwrite\gtoutfile
\long\gdef\makeheadfile{  
{\def\\{, }\def\s{ }
\immediate\openout\gtoutfile head.xxx
\immediate\write\gtoutfile{Proxy-for: \ifx\theasciiauthors\relax
\theauthors\else\theasciiauthors\fi\s<\ifx\theasciiemail\relax\theemail\else\theasciiemail\fi>}
\immediate\write\gtoutfile{\noexpand\\}
\immediate\write\gtoutfile{Authors: \ifx\theasciiauthors\relax
\theauthors\else\theasciiauthors\fi}
{\def\\{ }\immediate\write\gtoutfile{Title: \ifx\theasciititle\relax
\thetitle\else\theasciititle\fi}}
\immediate\write\gtoutfile{Subj-class: GT or SG, GR etc}
\immediate\write\gtoutfile{MSC-class: \theprimaryclass\ifx\thesecondaryclass\relax\else, \thesecondaryclass\fi}
\immediate\write\gtoutfile{Journal-ref: Algebr. Geom. Topol. \thevolumenumber\s
(\thevolumeyear) \startpage-\finishpage}
\immediate\write\gtoutfile{Comments: Published by Algebraic and
Geometric Topology at}
\immediate\write\gtoutfile{\s\s\s  http://www.maths.warwick.ac.uk/agt/AGTVol\thevolumenumber/agt-\thevolumenumber-\thepapernumber.abs.html}
\immediate\write\gtoutfile{\noexpand\\}
\immediate\write\gtoutfile{}
\ifx\theasciiabstract\relax
\immediate\write\gtoutfile{\theabstract}\else
\immediate\write\gtoutfile{\theasciiabstract}\fi
\immediate\write\gtoutfile{}
\immediate\write\gtoutfile{\noexpand\\}
\immediate\write\gtoutfile{}
\immediate\closeout\gtoutfile}}  

\def\maketitlepage{\makeagttitle\makeheadfile}
\let\makeshorttitle\maketitlepage
\let\maketitle\maketitlepage

\lognumber{56}
\volumenumber{5}
\volumeyear{2005}
\papernumber{56}
\pagenumbers{1433}{1450}
\received{20 July 2005} 
\accepted{2 September 2005}
\published{17 October 2005}
\usepackage{amsmath,amssymb}

\newtheorem{Theorem}{Theorem}
\newtheorem{Proposition}{Proposition}
\newtheorem{Lemma}{Lemma}
\newtheorem{Corollary}{Corollary}
\newtheorem{Claim}{Claim}
\theoremstyle{definition}
\newtheorem{Remark}{Remark}

\def\NN{{\mathcal N}}

\def\FF{{\mathcal F}}

\def\WW{{\mathcal W}}

\def\HH{{\mathcal H}}

\def\Z{{\mathbb Z}}

\begin{document}
\title{Degree one maps between small 3-manifolds\\and Heegaard genus}

\author{Michel Boileau\\Shicheng Wang}

\address{Laboratoire \'Emile Picard, CNRS UMR 5580, Universit\'e Paul 
Sabatier\\118 Route de Narbonne, F-31062 TOULOUSE Cedex 4, France}
\secondaddress{%
LAMA Department of Mathematics, Peking University\\Beijing 100871,
China}

\asciiaddress{Laboratoire Emile Picard, CNRS UMR 5580, Universite Paul 
Sabatier\\118 Route de Narbonne, F-31062 TOULOUSE Cedex 4, France\\and\\LAMA Department of Mathematics\\Peking University, Beijing 100871, China}

\asciiemail{boileau@picard.ups-tlse.fr, wangsc@math.pku.edu.cn}
\gtemail{\mailto{boileau@picard.ups-tlse.fr}, \mailto{wangsc@math.pku.edu.cn}}

\begin{abstract}
We prove a rigidity theorem for degree one maps between small
3-manifolds using Heegaard genus, and provide some applications and
connections to Heegaard genus and Dehn surgery problems.
\end{abstract}



\maketitle

\section{Introduction}

All terminology not defined in this paper is standard,
 see \cite{He} and \cite{Ja}.

 Let $M$ and $N$ be two closed, connected,
orientable 3-manifolds.  Let $H$ be a (not necessarily connected)
compact 3-submanifold of $N$. We say that a degree one map $f:
M\to N$ is a \emph{homeomorphism outside H} if $f: (M,
M-\text{int}f^{-1}(H), f^{-1}(H))\to (N, N-\text{int}H, H)$ is a
map between the triples such that the restriction $f|:
M-\text{int}f^{-1}(H) \to N-\text{int}H$ is a homeomorphism. We
say also that $f$ is a \emph{ pinch} and $N$ is obtained from $M$
by \emph{pinching $W=f^{-1}(H)$ onto $H$}.

Let $H$ be a compact 3-manifold (not necessarily connected). We
use $g(H)$ to denote the \emph{Heegaard genus} of $H$, that is the
minimal number of 1-handles used to build $H$.

We define \emph{$mg(H)={\text max}\{g(H_{i}),\, \text{$H_{i}$ runs
over components of $H$}\}.$} It is clear that $mg(H)\le g(H)$ and
$mg(H)= g(H)$ if $H$ is connected.

A path-connected subset $X$ of a connected $3$-manifold is said to carry
$\pi_{1}M$ if
the inclusion homomorphism $\pi_{1}X \to \pi_{1}M$ is surjective.

In this paper, any incompressible surface in a 3-manifold is
2-sided and is not the 2-sphere. A closed 3-manifold $M$ is {\it
small} if it is orientable, irreducible and if it contains no
incompressible surface.

It has been observed by Kneser, Haken and
Waldhausen (\cite{Ha}, \cite{Wa}, see
also \cite{RW} for a quick transversality argument) that a degree one
map $M \to N$ between two closed, orientable 3-manifolds is homotopic to a map
which is a homeomorphism outside a handlebody corresponding to one side of
a Heegaard splitting of $N$. This fact is known as ``any
degree one map between 3-manifolds is homotopic to a pinch".

A main result of this paper is the following rigidity theorem.

\begin{Theorem}\label{heegaard genus} Let $M$ and $N$ be two closed,
small 3-manifolds. If there is a degree one map $f: M\to N$ which is a
homeomorphism outside an irreducible submanifold $H\subset N$, then either:
\begin{enumerate}
\item There is a component $U$ of $H$ which carries $\pi_{1}N$
and such that $g(U) \geq g(N)$, or
\item $M$ and $N$ are homeomorphic.
\end{enumerate}
\end{Theorem}

\begin{Remark} Given $M$ and $N$  two non-homeomorphic
small 3-manifolds , Theorem \ref{heegaard genus} implies that $N$ cannot be
obtained from $M$ by a sequence
of pinchings onto submanifolds of genus smaller than $g(N)$.
However Theorem \ref{heegaard genus} does not hold when $M$
is not small. Below are easy examples:\begin{itemize}
\item Let $f: P\# N\to N$ be a degree one map defined
by pinching $P$ to a 3-ball in $N$. Then $f$ is a
homeomorphism outside the 3-ball, which is genus zero and does not carry
$\pi_{1}N$.

\item Let $k$ be a knot in a closed, orientable $3$-manifold $N$ and let
$F$ be a once punctured closed surface. Let
$M$ be the 3-manifold
obtained by gluing the boundaries of $F\times S^1$ and of  $E(k)$ in
such a way that $\partial F\times \{x \}$ is matched with the meridian of
$k$, $x\in S^1$. Then a degree one map $f:M\to N$ pinching $F\times S^1$
to a tubular neigborhood
$\NN(k)$ of $k$,  is a homeomorphism outside a handlebody of
genus 1. If $\pi_{1}N$ is not cyclic or tivial, then $g(\NN(k)) < g(N)$
and $\NN(k)$ does not carry $\pi_{1}N$.
\end{itemize}
\end{Remark}

The pinched part of a degree one map between closed, orientable non-homeo\-mor\-phic surfaces
 has incompressible boundary \cite{Ed}. The following
straigtforward
corollary of Theorem \ref{heegaard genus}
gives an analogous result for small 3-manifolds:

\begin{Corollary}\label{degree one map}
Let $M$ and $N$ be two  closed, small, non-homeomorphic
3-mani\-folds. Let $f: M\to N$ be a degree one map and let $V\cup
H=N$ be a minimal genus Heegaard splitting for $N$. Then the map
$f$ can be homotoped to be a homeomorphism outside $H$
such that $f^{-1}(H)$ is $\partial$-irreducible.
\end{Corollary}

\begin{Remark}
Corollary \ref {degree one map} remains true for any strongly irreducible
heegaard splitting of $N$. Then the
argument, using Casson-Gordon's result \cite{CG}, is  essentially the same
as \cite [Theorem 3.1]{Le}, even if in \cite {Le} it is only proved for the
case $M=S^3$ and $N$ a homotopy
3-sphere. The proof in \cite{Le} is based on his main result \cite [Theorem
1.3]{Le}, but one can also
use a direct argument from degree one maps.
\end{Remark}

Theorem \ref{heegaard genus} follows directly from  two  rather
technical Propositions (Proposition \ref{boundary irreducible} and
Proposition \ref{incompressible}). Theorem \ref{heegaard genus}
and its proof lead to some results about Heegaard genus of small
3-manifolds and Dehn surgery on null-homotopic knots.

\begin{Theorem}\label{small} Let $M$ be a closed, small
3-manifold.
Let $F \subset M $ be a closed, orientable surface (not necessary
connected) which cuts $M$ into finitely many compact, connected
3-manifolds $U_{1},\ldots, U_{n}$.
Then there is a component $U_i$  which carries $\pi_{1}M$ and
such that $g(U_i) \geq g(M)$.
\end{Theorem}

\begin{Remark}
In general (see \cite{La}) one has only the upper bound:
$$g(M) \leq  \sum_{i=1}^n g(U_i) ) - g(F).$$
\end{Remark}

Suppose that $k$ is a null-homotopic knot in a closed orientable
3-manifold $M$. Its unknotting number $u(k)$ is defined as the
minimal number of self-crossing changes needed to transform it
into a trivial knot contained in a 3-ball in $M$.

\begin{Theorem}\label{unknotting} Let $k$ be a null-homotopic knot in  a closed, small 3-manifold $M$. 
If $u(k)< g(M)$, then every
closed 3-manifolds obtained by a non-trivial Dehn surgery along $k$ is not small.
In particular $k$ is determined by its complements.
\end{Theorem}

This article  is organized as follows.

In Section \ref{boundincomp} we  state and prove Proposition
\ref{boundary irreducible} which is the first step in the proof of
Theorem \ref{heegaard genus}. The second step, given by
Proposition \ref{incompressible} is proved in Section
\ref{incompsurface}; then Theorem \ref{heegaard genus} follows
from these two propositions. Section \ref{heegaardsmall} is
devoted to the proof of Theorem \ref{small}, and Section
\ref{null-homotopic} to the proof of Theorem \ref{unknotting}.

\medskip
{\bf Acknowledgements}\qua We would like to thank both the referee and
Professor Scharlemann for their suggestions which enhance the
paper. The second author is partially supported by MSTC and NSFC.

\section{Making the preimage of $H$
$\partial$-irreducible}\label{boundincomp}

The first step of the proof of Theorem \ref{heegaard genus} is given by the
following proposition:

\begin{Proposition} \label{boundary irreducible}  Let $M$ and $N$ be
two closed,
connected, orientable, irreducible 3-manifolds which have the same
first Betti number, but are not homeomorphic.

Suppose there is a degree one map $f_0: M\to N$ which is a homeomorphism
outside a compact irreducible
3-submanifold $H_0\subset N$ with $\partial H_0\ne \emptyset$. Then
there is a degree one map $f: M\to N$ which is a homeomorphism
outside an irreducible submanifold $H \subset H_0$ such that:\begin{itemize}

\item $\partial H \not = \emptyset$;

\item $mg(H) \leq mg(H_0)$,

\item Any connected component of $f^{-1}(H)$ is either
$\partial$-irreducible or a 3-ball, and there is at least one
component of $f^{-1}(H)$ which is $\partial$-irreducible.
\end{itemize}
\end{Proposition}

\begin{Remark} Since $M$ is not homeomorphic to $N$ it is clear that at
least one component of
$f^{-1}(H)$ is not a 3-ball.
\end{Remark}

\proof In the whole proof, 3-manifolds $M$ and $N$ are supposed to
meet all hypotheses given in the first paragraph of Proposition
\ref{boundary irreducible}.

By the assumption
there is a degree one map $f_0: M\to N$ which is a
homeomorphism outside an irreducible submanifold  $H_0\subset N$
with $\partial H_0 \not = \emptyset$.

Let $\HH_{0}$ be the set of all 3-submanifolds $H\subset H_0$ such
that:\begin{itemize}

\item[(1)] There is a degree one map $f: M\to N$ which is a homeomorphism
outside $H$;

\item[(2)] $\partial H \not = \emptyset$;

\item[(3)] $mg(H)\leq mg(H_0)$;

\item[(4)] $H$ is irreducible.

\end{itemize}

For an element $H\in \HH_{0}$, its complexity is defined as a pair
$$c(H) = (\sigma(\partial H), \pi_0(H))$$
with the lexicographic order, and where $\sigma(\partial H)$ is
the sum of the squares of the genera of the components of
$\partial H$, and $\pi_0(H)$ is the number of components of $H$.

\medskip
{\bf Remark on $c(H)$}\qua The second term of $c(H)$ is not used in this section,
but will be used in the next two sections.
\medskip

Clearly  $\HH_{0}$ is not the empty set, since by  assumption $H_0 \in
\HH_{0}$.

A compressing disk for $\partial H$ in $H$ is a properly embedded $2$-disk $(D,
\partial D) \subset (H, \partial H)$
such that $\partial D = D \cap \partial H$ is an essential simple closed
curve on $\partial H$ (i.e. does not bound a
disk on $\partial H$). In the following
we shall denote by $H \backslash \NN(D)$ the compact 3-manifold obtained from
$H$ by removing an open product neighborhood of $D$. The operation of
removing such
neighborhood is called \emph{splitting} $H$ \emph{along}
$D$.

\begin{Lemma}\label{compressing disk} Let $H$ be a compact orientable
3-manifold and let $(D, \partial D) \subset (H, \partial H)$ be a compressing disk . Then
$mg(H_*)\le mg(H)$, where $H_* = H \backslash \NN(D) $ is obtained by
splitting $H$ along $D$.
Moreover $c(H_*)<c(H)$.
\end{Lemma}

\proof
By Haken's lemma for boundary-compressing disk (\cite{BO}, \cite{CG}), a minimal
genus Heegaard surface for $H$ can be isotoped to meet
$D$ along a single simple closed curve. It follows that $mg(H_*) \leq mg(H)$.

Since $\partial D$ is an essential simple closed curve on $\partial H$,
it is easy to see that $\sigma(\partial H_*)< \sigma(\partial H)$,
therefore $c(H_*) < c(H)$. \endproof

The proof of Proposition \ref{boundary irreducible} follows from
the following:

\begin{Lemma}\label{bdincomp} Let $H\in \HH_{0}$ be an element which
realizes the minimal complexity, then any component of
$f^{-1}(H)$ which is not a 3-ball is $\partial$-irreducible.
\end{Lemma}

\proof Let $W_{0} \subset W = f^{-1}(H)$ be a component
which is not homeomorphic to a 3-ball. Such a component exits
since $M$ is not homeomorphic to $N$.
To prove that $W_{0}$ is $\partial$-irreducible, we argue by
contradiction.

If $\partial W_{0}$ is compressible in $W$,
there is a compressing disc $(D, \partial D)\to (W,\partial W)$ whose
boundary is an essential simple closed curve on
$\partial W$.

Since $f: M\to N$ is a homeomorphism outside
the submanifold $H\subset N$ the restriction
$f|: (W,\partial W)\to (H, \partial H)$ maps $\partial W$
homeomorphically onto $\partial H$. Therefore
$f(\partial D)$ is an essential simple closed curve on $\partial H$  which
bounds the immersed disk $f(D)$ in $H$.
By Dehn's Lemma, $f(\partial D)$ bounds an embedded disc $D^*$ in $H$.

\begin{Lemma}\label{homotopy} By a homotopy of $f$, supported on $W =
f^{-1}(H)$ and
constant on $\partial W$, we can achieve that:
\begin{itemize}
\item  $f|: W \to H$  is a homeomorphism in a collar
neighborhood of $\partial W \cup D$,
\item $f|^{-1}(D^*)= D\cup S$, where $S$ is a closed orientable surface.
\end{itemize}
\end{Lemma}

\proof We define a homotopy $F: W\times [0,1] \to H$ by the following steps:

(1)\qua  $F(x,0)=f(x)$ for every $x\in W$;

(2)\qua  $F(x,t)=F(x,0)$ for every $x\in \partial f^{-1}(H) =
\partial W$ and for every $t\in [0,1]$;

(3)\qua Then we extend $F(x,1): D\times \{ 1\} \to D^*$  by a
homeomorphism.

We have defined $F$  on $D\times \{0\} \cup \partial D \times [0,1]\cup
D\times \{1\}$ which
is homeomorphic to a 2-sphere $S^2$. Since $H$ is irreducible, by the
Sphere theorem
$\pi_2(H) = \{0 \}$. Hence:

(4)\qua We can extend $F$ to $D\times [0,1]$;

Now $F$ has been defined on $W\times \{0\} \cup \partial W  \times
[0,1]\cup D\times [0,1]$,
which is a deformation retract of $W\times [0,1]$, therefore:

(5)\qua We can finally extend $F$ on $W\times [0,1]$.

After this homotopy we may assume that $f(x)=F(x, 1)$, for every
$x\in W$. Then by construction this new $f$ sends $\partial W \cup D$
homeomorphically to $ \partial H\cup D^*$.
By transversality, we may further assume that $f|: W\to H$ is a
homeomorphism in a
collar neighborhood of $\partial W  \cup D$ and that
$f|^{-1}(D^*)= D\cup S$, where $S$ is a closed surface. \endproof

The following lemma will be useful:

\begin{Lemma}\label{homology} Suppose $f: M\to N$ is a degree one map
between two closed orientable 3-manifolds with the same first Betti number
$\beta_1(M)=\beta_1(N)$.
Then $f_{\star} : H_{2}(M;\Z) \to H_{2}(N;\Z)$ is an isomorphism.
\end{Lemma}

\proof Since $f: M\to N$ is a degree one map, by \cite [Theorem
I.2.5]{Br}, there is a homomorphism $\mu: H_2(N;\Z)\to H_2(M;\Z)$
such that $f_{\star}\circ \mu : H_2(N;\Z)\to H_2(N;\Z)$ is the
identity, where $f_{\star} : H_2(M;\Z)\to H_2(N;\Z)$ is the
homomorphism induced by $f$.

In particular $f_{\star} : H_2(M;\Z)\to H_2(N;\Z)$ is surjective. Then the injectivity  follows from the fact that
$H_2 (M;\Z)$ and $H_2(N;\Z)$ are torsion free abelian groups with the same
finite rank $\beta_2(M)=\beta_1(M)=\beta_1(N)=\beta_2(M)$.\endproof

Since the degree one map $f : M \to N$ is a homeomorphism outside $H$, 
the Mayer-Vietoris sequence and
Lemma \ref{homology} imply that  $f_{\star} : H_{2}(W;\Z) \to H_{2}(H;\Z)$ is an
isomorphism.

Let $S'$ be a connected component of $S$. Since $f(S')\subset D^*$, the
homology class $[f(S')]=f_{\star}([S'])$ is zero in $H_2(H,\Z)$. Hence
the homology class $[S']$ is zero in $H_2(W,\Z)$, because
$f_{\star}:H_2(W,\Z)\to H_2(H,\Z)$ is an isomorphism.
It follows that  $S'$ is the boundary of a compact submanifold of $W$.
Therefore $S'$ divides $W$ into two parts
$W_1$ and $W_2$ such that $\partial W_2 = S'$ and $W_1$ contains $\partial
W\cup D$.

We can define a map $g: W\to H$ such that:

(a)\qua $g|_{W_1}=f|_{W_1}$ and $g(W_2)\subset D^*$.

Then by slightly pushing the image $g(W_2)$ to the correct side of $D^*$,
we can improve the map  $g: W\to H$ such that:

(b)\qua $g|\partial W = f|\partial W $,

(c)\qua $g^{-1}(D^*)=D \cup (S \setminus S')$ and $g: \NN(D)\to
\NN(D^*)$ is a homeomorphism.

After finitely many such steps we get a map $h: W\to H$ such that:

(b)\qua $h|\partial W = f|\partial W $,

(d)\qua $h^{-1}(D^*)=D$ and $h: \NN(D)\to \NN(D^*)$ is a homeomorphism.

Let $H_* = H \backslash \NN(D)$ obtained by splitting $H$ along $D$.
Then $H_*$ is still an irreducible 3-submanifold of
$N$ with $\partial H_* \not = \emptyset$.

Now $f|_{M-\text{int}W}$ and $h|_{W}$ together provide a degree
one map $h: M\to N$. The transformation from
$f$ to $h$ is supported in $W$, hence $h$ is a homeomorphism
outside the irreducible submanifold  $H_*$ of $N$.

Since $H_*$ is  obtained by splitting $H$ along a compressing disk,
we have $H_* \subset H_0$ and $H_*$ belongs to $\HH_{0}$. Moreover $mg(H_*)
\leq mg(H)$ and $c(H_*)< c(H)$ by
Lemma \ref{compressing disk}.

This contradiction finishes the proof of Lemma \ref{bdincomp} and
thus the proof of Proposition \ref{boundary irreducible}.
\endproof

\section{Finding a closed incompressible surface in the
domain}\label{incompsurface}

Since closed, orientable, small 3-manifolds are irreducible
and have first Betti number equal to zero, Theorem \ref{heegaard genus} is
a direct corollary of the following proposition:

\begin{Proposition}\label{incompressible} Let $M$ and $N$ be two closed,
connected, orientable, irreducible 3-manifolds whith the same
first Betti number.
Suppose that there is a degree one map $f: M\to N$ which is a homeomorphism
outside an irreducible
submanifold $H_0\subset N$  such that for
each connected component $U$ of $H_0$, either $g(U) < g(N)$ or $U$ does
not carry $\pi_{1}N$. Then either
$M$ contains an incompressible orientable surface or $M$ is homeomorphic to
$N$.
\end{Proposition}

Let $(M,N)$ be a pair of closed orientable 3-manifolds such that there is a
degree one map from $M$ to $N$. We say that condition $(*)$ holds for
the pair $(M,N)$ if:
$$ (*)  \qquad \pi_{1}N = \{1\} \, \, \, \, \, \text{implies} \, \, \, \,
M  =
S^3.$$
For the proof we first assume that condition $(*)$ holds for the pair
$(M,N)$.

\medskip

{\bf Proof of Proposition \ref{incompressible} under condition $(*)$}

By the assumptions, there is a degree one map
$f: M\to N$ which is a homeomorphism outside
an irreducible submanifold $H_{0}\subset N$ with $\partial H_{0} \not =
\emptyset$ and such that for
each connected component $U$ of $H_0$ either $g(U) < g(N)$ or $U$ does not
carry $\pi_{1}N$. We assume
moreover that $M$ is not homeomorphic to $N$. Our goal is to show that $M$
contains an incompressible surface.

Similar to Section \ref{boundincomp}, let $\HH$ be the set of all
3-submanifolds $H \subset N$  such that: \begin{itemize}
\item[(1)] There is a degree one map $f: M\to N$ which is a homeomorphism
outside $H$.
\item[(2)] $\partial H$ is not empty.
\item[(3)] For each component $U$ of $H$,
either $g(U) < g(N)$ or $U$ does not carry $\pi_{1}N$.
\item[(4)] $H$ is irreducible.
\end{itemize}

The set $\HH$ is not empty  by our assumptions.

The complexity $c(H) = (\sigma(\partial H), \pi_0(H))$ for the elements of
$\HH$ is defined like in Section \ref{boundincomp}.

\begin{Lemma}\label{ball} Assume that there is a degree one map $f: M\to N$
which is a homeomorphism outside a submanifold $H\subset N$. If
$H$ contains 3-ball component $B^3$, then $f$ can be homotoped to
be a homeomorphism outside $H_*$, where $H_*= H-B^3$. Moreover if
$H$ is irreducible, then $H_*$ is also irreducible.
\end{Lemma}

\proof By our assumption, there is a degree one map $f: M\to N$
which is a homeomorphism outside a submanifold $H\subset N$ and
$H$ contains a $B^3$ component. Since $f|: f^{-1}(\partial H) \to
\partial H$ is a homeomorphism, then $f^{-1}(\partial B^3)$ is a
2-sphere $S^2_* \subset M$. Since $M$ is irreducible, $S^2_*$
bounds a 3-ball $B_*^3$ in $M$. Then either \begin{itemize}
\item[(a)] $M-\text{int}f^{-1}(B^3)=B_*^3$, or
\item[(b)] $f^{-1}(B^3)=B^3_*$.
\end{itemize}

In case (a), $N= f(B_*^3)\cup B^3$ is a union of two homotopy
3-balls with their boundaries identified homeomorphically, and
clearly $\pi_{1}N = \{1\}$.  So $M= S^3$ by assumption $(*)$. Hence
(b) holds in either case.

In case (b), by a homotopy of $f$ supported in $f^{-1}(B^3)$, we
can achieve  that $f|: f^{-1}(B^3)\to B^3$ is a homeomorphism.
Then $f$ becomes a homeomorphism outside the irreducible
3-submanifold $H_* \subset N$, obtained from $H$ by deleting the
3-ball $B^{3}$.

The last sentence in Lemma \ref{ball} is obviously true.
\endproof

Let $H\in \HH$ be an element which realizes the minimal
complexity. By Lemma \ref{ball} no component of $H$ is a 3-ball,
hence no component of $\partial H$ is a 2-sphere since $H$ is
irreducible. Therefore no component of $f^{-1}(H)$ is a 3-ball and
$\partial f^{-1}(H)$ is incompressible in $f^{-1}(H)$ by the proof
of Lemma \ref{bdincomp}.

Since $f : M-\text{int}f^{-1}(H) \to  N-\text{int}H$ is a homeomorphism,
$\partial f^{-1}(H)$ is incompressible in
$M-\text{int}f^{-1}(H)$ if and only if  $\partial H$ is
incompressible in $N-\text{int}H$. For simplicity we will set $V =
N-\text{int}H$, then $N=V\cup H$.

Then the proof of Proposition \ref{incompressible} under condition
$(*)$ follows from:

\begin{Lemma}\label{minimal}  If $\partial H$ is compressible in $V$, then
there is $H_*\in \HH$ such that $c(H_*) < c(H)$.
\end{Lemma}

\proof  Suppose $\partial H$ is compressible in $V$. Let $(D,\partial
D)\subset (V,\partial V)$ be a compressing disc.
By surgery along $D$, we get two submanifolds  $H_1$ and $V_1$ as follows:
$$H_1=H\cup \NN(D),\qquad V_1=V \backslash \NN(D).$$
Since $H_1$ is obtained from $H$ by adding a 2-handle, for each
component $U'$ of $H_1$ there is a  component $U$ of $H_0$ such
that $g(U') \leq g(U)$ and $\pi_{1}U'$ is a quotient of
$\pi_{1}U$, hence $H_1$ verifies the defining condition (3) of
$\HH$. Moreover $f$ is still a homeomorphism outside $H_{1}$
because $H_1$ contains $H$ as a subset. Clearly $\partial H_{1}
\not = \emptyset$. Hence $H_1$ satisfies also the defining
conditions (1) and (2)  of $\HH$. We notice that $c(H_1) < c(H)$
because $\sigma(\partial H_1) < \sigma(\partial H)$.

We will modify $H_1$ to become  $H_*\in \HH$ with $c(H_*) \leq c(H_1)$.
The modification will be divided into two steps carried by Lemma \ref{sphere}
and Lemma \ref{irreducible}  below. First the following standard lemma will
be useful:

\begin{Lemma}\label{standard} Suppose $U$ is a connected 3-submanifold in
$N$ and let  $B^3\subset N$ be a 3-ball with $\partial B^3 = S^2$.

{\rm(i)}\qua Suppose $S^2 \subset \partial U$. If $\text{int}\, U\cap
B^3\ne \emptyset$, then
$U\subset B^3$. Otherwise  $U\cap B^3=S^2$.

{\rm(ii)}\qua if $\partial U\subset B^3$, then either $U\subset B^3$, or
$N-\text{int}\, U\subset B^3$.
\end{Lemma}

\proof For (i): Suppose first $\text{int}\, U\cap B^3\ne \emptyset$.
Let $x\in \text{int}\, U \cap B^3$. Since $U$ is connected,
then for any $y \in U$, there is a path $\alpha\subset U$ connecting
$x$ and $y$. Since $S^2$ is a component of $\partial U$, $\alpha $
does not cross $S^2$. Hence $\alpha\subset B^3$ and $y\in B^3$,
therefore $U\subset B^3$.

Now suppose  $\text{int}\, U\cap B^3= \emptyset$. Let $x\in \partial U\cap
B^3$. If $x\in \text{int}\, B^3$,
then there is $y\in \text{int}\, U\cap B^3$. It contradicts the assumption.
So $x\in \partial B^3=S^2$.

For (ii): Suppose that  $U$ is not a subset of $B^3$, then there
is a point $x\in U\cap (N-\text{int}B^3)$. Let $y\in N-\text{int}\, U$. If
$y\in N-\text{int}\, B^3$, there is a path
$\alpha$ in $N-\text{int}\, B^3$ connecting $x$ and $y$, since
$N-\text{int}\, B^3$ is connected. This path $\alpha$ does not meet
$\partial U$, because  $\partial U\subset B^3$. This would
contradict that $x\in U$ and $y\in N-\text{int}\, U$. Hence we must
have $y\in B^3$, and therefore $N-\text{int}\, U\subset B^3$. \endproof

\begin{Lemma}\label{sphere} Suppose $H_{1}$ meets the defining conditions
(1), (2) and (3) of the set $\HH$. Then $H_1$ can be modified
to be a 3-submanifold $H_*\subset N$ such that:

{\rm(i)}\qua $\partial H_*$ contains no 2-sphere;

{\rm(ii)}\qua $c(H_*) \leq c(H_1)$;

{\rm(iii)}\qua $H_*$ still meets the the defining conditions (1) (2) (3)
of $\HH$.
\end{Lemma}

\proof We suppose that $\partial H_1$ contains a 2-sphere
component $S^2$, otherwise Lemma \ref{sphere} is proved. Then $S^2$ bounds
a 3-ball $B^3$ in $N$ since $N$ is irreducible. We consider two cases:

\medskip
{\bf Case (a)}\qua \textsl{$S^2$ does bound a 3-ball $B^3$ in $H_1$.}

In this case $B^3$ is a component of $H_1$. By Lemma \ref{ball}, $f$ can
be homotoped to be a homeomorphism outside $H_{2}= H_{1}-B^3$.

\medskip
{\bf Case (b)}\qua \textsl{$S^2$ does not bound a 3-ball $B^3$ in $H_1$.}

Let $H'_1$ be the  component of $H_1$ such that $S^2 \subset
\partial H'_1$.
By Lemma \ref{standard} (i), either:

{\bf (b$'$)}\qua  $H'_1 \subset B^3$, or

{\bf (b$''$)}\qua  $H'_1 \cap B^3 = S^2$.

\medskip {\bf In case (b$'$)}, let $H_2 = H_1 - B^3$.  By Lemma
\ref{ball} $f$ can be homotoped to be a homeomorphism outside
$H_2$.  Note $H_2\ne \emptyset$, otherwise  $M$ and $N$ are
homeomorphic, which contradicts our assumption.

\medskip{\bf In case (b$''$)}, let $H_{2}=H_{1}\cup B^3$, then $\partial
H_{2}$ has at least one
component less than $\partial H_1$. Since we are
enlarging $H_1$, $f$ is a homeomorphism outside $H_{2}$.

It is easy to check that in each case (a), (b$'$), (b$''$) the
components of $H_2$ verify the defining condition (3) of $\HH$ and
$c(H_2) \le c(H_1) < c(H)$. Moreover $H_2$ is not empty because
$M$ and $N$ are not homeomorphic, and  $\partial H_2 \not =
\emptyset$ since $g(H_2)\leq g(H_1)< g(N)$. Hence each of the
transformations (a), (b$'$) and (b$''$) preserves properties (ii) and
(iii) in the conclusion of Lemma \ref{sphere}. Since  each one
strictly reduces the number of components of $H_1$ or of $\partial
H_1$, after a finite
 number of such transformations we  reach a 3-submanifold $H_*$ of $N$
such that $H_*$ meets the properties (ii) and (iii) of Lemma \ref{sphere},
and $\partial H_*$ contains no 2-sphere components. This proves Lemma
\ref{sphere}.\endproof

\begin{Lemma}\label{irreducible}  Suppose that $H_{1}$ meets conditions (i),
(ii) and (iii) in the conclusion of Lemma \ref{sphere}. Then $H_1$ can be
modified to be a 3-submanifold $H_*$ of $N$ such that:

{\rm(a)}\qua $H_*$ is irreducible;

{\rm(b)}\qua $c(H_*) \leq c(H_1)$ is not increasing;

{\rm(c)}\qua $H_*$ still meets the the defining conditions (1), (2), (3)
of $\HH$.

In particular $H_*$ belongs to $\HH$.
\end{Lemma}

\proof
If there is an essential 2-sphere $S^2$ in $H_1$, it must separate $N$
since $N$ is irreducible. Let $H'_1$ be the
component of $H_1$ containing $S^2$. The 2-sphere $S^2$ induces a connected
sum decomposition of $H'_{1}$: it separates
$H'_{1}$ into two connected parts
$K_\circ$ and $K'_\circ$, such that:
$$H'_{1} = K\#_{S^2}K' = K_{\circ}\cup _{S^2}K'_{\circ},$$
\noindent $K_{\circ}\subset H_1$ (resp.
$K'_{\circ}\subset H_1$) is homeomorphic to a once punctured
$K$  (resp. a once punctured $K'$).

By Haken's Lemma, we have:
 $$g(H'_{1})=g(K) +g(K').$$
Neither $K_{\circ}$ nor $K'_{\circ}$ is a $n$-punctured 3-sphere, $n\ge 0$,
because $\partial H_1$
contains no 2-sphere component, hence:
 $$g(K)< g(H'_{1}) \,\qquad \text{and} \,\qquad g(K') < g(H'_{1})$$
Since $N$ is irreducible, $S^2$ bounds a 3-ball $B^3$ in $N$. We may assume
that
$\text{int} K_{\circ}\cap B^3 = \emptyset$ and
$\text{int} K'_{\circ}\cap B^3\ne \emptyset$. By Lemma \ref{standard} (i),
we have $K_{\circ}\cap B^3=S^2$ and
$K'_{\circ}\subset B^3$.

Moreover $\partial H'_{1}\cap B^3\ne\emptyset$, otherwise $K'_{\circ}$ is
homeomorphic to $B^3$,
in contradiction with the assumption
that $S^2$ is a  2-sphere of connected sum.

\begin{Lemma}\label{essential} $\partial H'_{1}$ is not a subset of $B^3$.
\end{Lemma}

\proof
We argue by contradiction. If $\partial H'_{1}$ is a subset of $B^3$, we
have $N-\text{int}H'_{1}\subset B^3$ by Lemma \ref{standard} (ii), since
$H'_{1}$ is not a subset of $B^3$. Then:
 $$N = H'_{1}\cup (N-\text{int}H'_{1}) = H'_{1}\cup B^3
=(K_\circ\#_{S^2}K'_\circ)\cup B^3 = K_\circ \cup_{S^2}B^3=K.$$
Hence $K$ is homeomorphic to the whole $N$. If $g(H'_1) <g(N)$, this
contradicts the fact that $g(K) < g(H'_1) <g(N)$.
If $H'_1$ does not carry $\pi_{1}N$ this contadicts the fact that $K
\subset H'_1$. \endproof

By Lemma \ref{essential}, $\partial H'_{1}$ (and therefore $\partial H_1$)
has components  disjoint from $B^3$. Therefore if we replace $H_1$ by
$H_2=H_1\cup B^3$, then $\partial H_2$ is  not empty and it has no
component which is a 2-sphere. Moreover the application of Haken's
Lemma above shows that $g(H_2) < g(H_1)$.

Since we are enlarging $H_1$,  $f$ is  a homeomorphism outside $H_2$, and
clearly $H_2$ still meets the the defining condition
(3) of $\HH$. Moreover $c(H_2) \leq c(H_1)$. Hence the transformation from
$H_1$ to $H_2$ preserves properties (b) and (c) in
the conclusion of Lemma \ref{irreducible}.
Since it strictly reduces $g(H_1)$, after a finite
number of such transformations we will reach a 3-submanifolds $H_* \subset N$
such that $H_*$ meets conditions (b) and (c)
in the conclusion of Lemma \ref{irreducible}, but does not contain any
essential 2-sphere. This proves Lemma \ref{irreducible}.\endproof

Lemma \ref{sphere} and Lemma \ref{irreducible} imply Lemma \ref{minimal}.
Hence we have proved Proposition \ref{incompressible}  under condition
$(*)$. \endproof

{\bf Proof of Proposition \ref{incompressible}}\qua Let $M$
and $N$ be two  closed, small 3-manifolds which are not
homeomorphic. Suppose there is  degree one map $f: M\to N$ which
is a homeomorphism outside an irreducible submanifold $H\subset N$
such that: for each  component $U$ of $H$, either $g(U) < g(N)$ or
$U$ does not carry $\pi_{1}N$.

Condition $(*)$ in the above proof of  Proposition
\ref{incompressible} is only used in
the proof of Lemma \ref{ball}, when $H$ contains a 3-ball component
$B^3$ and that 
$M-\text{int}f^{-1}(B^3)=B_*^3$ and
$f^{-1}(B^3) \not = B^3_*$. Indeed we can now prove that this case
cannot occur.

If this  case happens then $\pi_{1}N = \{1\}$ and thus $mg(H) <
g(N)$, since every  component of $H$ carries $\pi_{1}N$. By
replacing $f^{-1}(B^3)$ by a 3-ball $B^3_\#$,  we obtain a
degree one map $\bar f: S^3=B^3_*\cup B^3_\#\to N$ defined by
$\bar f| B_* = f |B_*$ and $\bar f|: B^3_\#\to B_3$ is a
homeomorphism.
Then $\bar f: S^3\to N$ is a map which is a homeomorphism outside
a submanifold $H'=H-B^3$. Clearly $mg(H')= mg(H) < g(N)$. Furthermore condition
$(*)$ now holds.

Since Proposition \ref{incompressible} has been proved under condition
$(*)$, we have that $N$ must be homeomorphic to $S^3$, since $S^3$ does not
contain any incompressible surface. It
would follow
that $mg(H)<0$, which is impossible.

The proof of Proposition \ref{incompressible}, and hence of Theorem
\ref{heegaard genus} is now complete.\endproof

\section{Heegaard genus of small 3-manifolds} \label{heegaardsmall}

This section is devoted to the proof of Theorem \ref{small}.

Let $M$ be a closed orientable irreducible 3-manifold. Let $F
\subset M $ be a closed orientable surface (not necessary
connected) which splits $M$ into finitely many compact connected
3-manifolds $U_{1},\ldots, U_{n}$.

Let $M \backslash \NN(F)$ be the manifold $M$ split along the surface $F$. We
define the complexity of the pair $(M, F)$ as
$$c(M,F) = \{\sigma(F), \pi_0(M \backslash \NN(F))\},$$
where $\sigma(F)$ is the sum of the squares of the genera of the
components of $F$ and $\pi_0(M \backslash \NN(F))$ is the number
of components of $M \backslash \NN(F)$.

Let $\FF$ be the set of all closed surfaces $F$ such that for each
component $U_i$ of $M \setminus F$, either $g(U_i) < g(M)$ or
$U_i$ does not carry $\pi_{1}M$.

\begin{Remark} This condition implies that the surface $F \not = \emptyset$
for every $F \in \FF$.
\end{Remark}

With the hypothesis of Theorem \ref{small}, the set $\FF$ is not
empty. Let $F\in \FF$ be a surface realizing the minimal
complexity. Then the following Lemma implies Theorem \ref{small}.

\begin{Lemma}\label{nosphere}
A surface $F\in \FF$ realizing the minimal complexity  contains no 2-sphere
component and is incompressible.
\end{Lemma}

\proof
The arguments are analogous to those used in the proof of Propositions
\ref{incompressible}. We argue by contradiction.

Suppose that $F$ contains a 2-sphere component $S^2$. It bounds a
3-ball $B^3 \subset M$, since $M$ is irreducible. Let $U_1$ and
$U_2$ be the closures of the components of $M\setminus \NN(F)$
which contain $S^2$. Then by Lemma \ref{standard} (i),
either:\begin{itemize}
\item $U_2\subset B^3$ and $U_1 \cap B^3 = S^2$, or
\item $U_1\subset B^3$ and $U_2 \cap B^3 = S^2$.
\end{itemize}

Since those two cases are symmetric, we may assume that we are in the first
case.
We consider the surface $F'$ corresponding to the decomposition
$\{U'_1,\ldots,U'_k \}$ of $M$
with  $U'_1 = U_1 \cup B^3$, after forgetting all $U_i \subset B^3$ and
then re-indexing the remaining $U_i$'s to be $U'_2,\ldots, U'_k$.
This operation does not increase the Heegaard genus of any one of the
components of the new decomposition.
Moreover if $U_1$ does not carry $\pi_{1}M$, the same holds for $U'_{1}$.
Hence $F'$ still belongs to $\FF$.
However, this operation strictly decreases the number of components of $F$,
hence $c(F')<c(F)$, in contradiction with our choice of $F$.

Suppose that the surface $F$ is compressible. Then some essential simple
closed curve $\gamma$ on $F$ bounds an embedded disk in $M$. Let $D'$ be a
such a
compression disk with the minimum number of circles of intersection in
$\text{int} D'\cap F$.
Then a subdisk of $D'$ bounded
by an innermost circle of intersection is contained inside one of the
$U_i$, say $U_1$.

Let  $(D, \partial D)\subset (U_1, F\cap \partial U_{1})$ be such an
innermost disk. Let $U_2$ be adjacent to $U_1$ along $F$, such that
$\partial{D} \subset \partial U_{2}$.
By surgery along $D$, we get a new surface $F'$ which gives a new
decomposition $\{U_1',\ldots,U_n'\}$ of $M$ as follows:
$$U'_1=U_1 \backslash \NN(D) ,\qquad U'_2= U_2\cup \NN(D), \qquad U'_i=U_i,
\, \text{for} \, i \geq 3.$$
Then  $g(U'_i) \leq g(U_i)$, for $i = 1,\dots, n$. Moreover if $U_i$ does
not carry
$\pi_{1}M$, the same holds for $U'_{i}$. Hence $F' \in \FF$.
However, $\sigma(F') < \sigma(F)$ since $\partial D$ is an essential circle
on $F$. Therefore $c(F') < c(F)$ and we reach a contradiction.\endproof

\section{Null-homotopic knot with small unknotting
number}\label{null-homotopic}

In this section we prove Theorem \ref{unknotting}.

Suppose $M$ is a closed, small 3-manifold and
$k\subset M$ is a null-homotopic knot with $u(k)< g(M)$. Then clearly $M$
is not the 3-sphere.

If $k$ is a non-trivial knot in a 3-ball $B^3 \subset M$. Then the compact
3-manifold $B^3(k,\lambda)$ obtained by any
non-trivial surgery of slope
$\lambda$ on $k$ will not be a 3-ball by \cite{GL}. Therefore
$M(k,\lambda)$ contains an essential 2-sphere.

Hence below we assume that $k$ is not contained in a 3-ball.

Since the knot $k\subset M$ is null-homotopic with unknotting number
$u(k)$, $k$ can be obtained
from a trivial knot $k' \subset B^3 \subset M$ by
$u(k)$ self-crossing changes. Let $D' \subset M$ be an embedded disk
bounded by $k'$.
If we let $D'$ move following the
self-crossing changes from $k'$ to $k$, then each self-crossing change
corresponds
to an identification of pairs of arcs in $D'$.
Hence one obtains a singular disk $\Delta$ in $M$ with $\partial \Delta =
k$ and
with $u(k)$ clasp singularities. Since $\Delta$
has the homotopy type of a graph, its regular neighborhood $\NN(\Delta)$
is a handlebody of genus $g(\NN(\Delta)) = u(k) < g(M)$.

First we prove the following lemma  which is a particular case of a more
general result about Dehn surgeries on
null-homotopic knots, obtained in \cite{BBDM}. Since this paper is not yet
available, we give
here a simpler proof in this particular case.

\begin{Lemma}\label{nonhomeo} With the hypothesis above, if the slope
$\alpha$ is not the meridian slope of $k$, then $M(k,\alpha)$ is not
homeomorphic to $M$.
\end{Lemma}

\proof
Since $M$ is irreducible and
$k\subset M$ is not contained in a 3-ball,
$M-\text{int}\NN(k)$ is irreducible and $\partial$-irreducible. Hence $1
\leq u(k) < g(M)$ and $M$ cannot be a lens space.

Let consider the set $\WW$ of compact, connected, orientable,
3-submanifolds $W \subset M$ such
that:\begin{enumerate}
\item $k \subset W$ is null-homotopic in $W$;
\item there is no 2-sphere component in $\partial W$;
\item $g(W) < g (M)$.
\end{enumerate}

By hypothesis the set $\WW$ is not empty since a regular neighborhood
$\NN(\Delta)$ of a singular unknotting disk for
$k$ is a handlebody of genus $\geq 1$.

\begin{Claim} For a 3-submanifold $W_0 \in \WW$ with a minimal complexity
$c(W_0) = \sigma(\partial W_0)$, the surface
$\partial W_0$ is incompressible in the exterior $M-\text{int}\NN(k)$.
\end{Claim}

\proof
If $\partial W _0$ is compressible in $M -\text{int} W _0$, let
$(D, \partial D) \hookrightarrow (M -\text{int} W _0, \partial W_0)$
be a compression disk for $\partial W _0$.
The $3$-manifold $W_1 = W_0 \cup \NN(D)$, obtained by adding a 2-andle to
$W_0$, is a compact, connected submanifold
of $M$ containing $k$.

Any
2-sphere in
$\partial W_1$  bounds a 3-ball in  $M-\text{int}\NN(k)$ since it is
irreductible. Hence after gluing some 3-ball along the
boundary, we may  assume that $W_1$ contains no 2-sphere component.
Moreover $k \subset W_1$ is null-homotopic in $W_1$ and
$g(W_1) \leq g(W_0) < g (M)$. It follows that $W_1 \in \WW$. Since $c(W_1)
< c(W_0)$ we get a
contradiction.

If $\partial W _0$ is compressible in $W_0 -\text{int}\NN(k)$, let
$(D, \partial D) \hookrightarrow , (W_0 -\text{int}\NN(k),
\partial W_0)$ be a compression disk for $\partial W _0$. Let
$W_2$ be the component of the 3-manifold $W_0 \backslash \NN(D)$
which contains $k$. As above, after possibly gluing some 3-ball
along the boundary, we may assume that $\partial W_2$  contains no
2-sphere component. The knot $k \subset W_2$ is null-homotopic in
$W_2$, since it is null-homotopic in $W_0$ and $\pi_{1}W_2$ is a
factor of the free product decomposition of $W_0$ induced by the
$\partial-compression$ disk $D$. Moreover by Lemma
\ref{compressing disk} $g(W_2) \leq g (W_0) < g (M)$. It follows
that $W_2 \in \WW$ and $c(W_2) < c(W_0)$. As above this
contradicts the minimality of $c(W_0)$.\endproof

To finish the proof of Lemma \ref{nonhomeo} we distinguish two cases:

\medskip
{\bf (a)}\qua {\sl The surface  $\partial W_0$ is compressible in
$W_0(k,\alpha)$}\qua Then one can apply Scharlemann's
theorem \cite [Thm 6.1]{Sch}. The fact that $k \subset W_0$ is
null-homotopic rules out cases a) and b) of Scharlemann's theorem.
Moreover by \cite [Prop.3.2]{BW} there is
a degree one map $g : W_0(k,\alpha) \to W_0$, and thus there is a simple closed
curve on $\partial W_0$ which is a  compression curve  both in
$W_0(k,\alpha)$  and in $W_0$. Therefore case d) of Scharlemann's theorem
cannot occure. The remaining case c) of Scharlemann's theorem shows that $k
\subset W_0$ is a non-trivial cable of a knot
$k_0 \subset W_0$ and that the surgery slope $\alpha$ corresponds to the
slope of the cabling annulus. But then the manifold
$M(k,\alpha)$ is the connected sum of a non-trivial Lens space with a
manifold obtained by Dehn surgery
along $k_0$. If $M(k,\alpha)$ is homeomorphic to the small 3-manifold $M$,
then $M$ and $M(k,\alpha)$ both
would be homeomorphic to a Lens space, which is impossible since $1\leq
u(k) < g(M)$.

\medskip
{\bf (b)}\qua {\sl The surface  $\partial W_0$ is incompressible in
$W_0(k,\alpha)$}\qua Since  $\partial W_0$ is incompressible in
$M-\NN(k)$, it is incompressible in $M(k,\alpha)$. Therefore $M(k,\alpha)$
and $M$ cannot be homeomorphic since $M$ is a small
manifold. \endproof

It follows from \cite [Prop.3.2]{BW} that there is a degree one
map $f : M(k,\alpha) \to M$ which is a homeomorphism outside
$\NN(\Delta)$. Since $g(\NN(\Delta)) = u(k) < g(M)$, Theorem
\ref{unknotting} is a consequence of Theorem \ref{heegaard genus}
and Lemma \ref{nonhomeo}. \endproof

\Addresses\recd

\end{document}